\documentclass[12pt]{amsart}
\usepackage{amssymb,amsmath,amsthm,latexsym}
\usepackage{mathrsfs}
\usepackage{mdwlist}
\usepackage{a4wide}

\theoremstyle{definition}

\newtheorem*{theorem*}{{\bf Main Theorem}}
\newtheorem*{remark*}{Remark}

\theoremstyle{remark}

\newcounter{smallromans}

  {\end{list}}

\newcounter{smallalphs}

  {\end{list}}

\author{Tomasz Kania}
\address{Mathematics Institute, University of Warwick, Gibbet Hill Rd, Coventry, CV4 7AL, England}
\email{tomasz.marcin.kania@gmail.com}
\title{On bases that are closed under multiplication}
\dedicatory{To appear in the American Mathematical Monthly}
\thanks{The author acknowledges with thanks funding received from the European Research Council / ERC Grant Agreement No.~291497.}
\begin{document}
\begin{abstract}It is well known that there is no basis of the field for real numbers regarded as a vector space over any proper subfield that is closed under multiplication. Mabry has extended this result to bases of arbitrary proper field extensions. The aim of this short communication is to notice that the proof of the result concerning the reals may be adjusted to a larger class of algebras (including full matrix algebras); thereby we subsume Mabry's result. 
\end{abstract}
\maketitle
Let $k$ be a field and let $A$ be a $k$-algebra, \emph{i.e.}, $A$ is a vector space over the field $k$ which is at the same time a (possibly non-unital or non-commutative) ring whose operations are compatible with the vector-space operations. (A ring is \emph{unital} when it is has an identity element with respect to multiplication.) For example, for every fixed natural number $n$, the Cartesian product $k^n$ of $n$ copies of $k$ is a paradigm example of a unital $k$-algebra when furnished with the operations of addition, multiplication and scalar multiplication defined component-wise.\smallskip

Since $A$ is a vector space over $k$ on its own, we may talk about bases in $A$ (that is, maximal linearly independent subsets of $A$). By linear algebra, a subset $H\subset A$ is a~basis if and only if every vector $x\in A$ may be written uniquely as a linear combination of vectors in $H$. (By the Kuratowski--Zorn lemma, every vector space has a basis.) Let us say that a basis $H$ of $A$ is \emph{closed under multiplication} if $x\cdot y\in H$ whenever $x$ and $y$ are in $H$. \smallskip 

Such bases have been studied by Mabry (\cite{mabry}) in the context of field extensions. More specifically, Mabry proved that if $k_1$ is a proper subfield of another field $k_2$, then no basis of the latter field, regarded as a $k_1$-algebra, is closed under multiplication. L\'{o}pez-Permouth, Moore and Szabo (\cite{lopez}) studied bases in unital algebras that consist of invertible elements and they have obtained an alternative proof of Mabry's result on field extensions. The aim of this note is to extend Mabry's result from a~different angle---we shall give a quick proof of a~more general yet still very simple result that will encompass not only field extensions but many other algebras too. Before we do so, we introduce a piece of terminology.\smallskip

An \emph{ideal} of a $k$-algebra $A$ is a linear subspace of $A$ that is closed under multiplying its elements by arbitrary elements of $A$ from left and right. The $k$-codimension of a~subspace $V$ of a vector space $W$ over $k$ is the $k$-dimension of the quotient space $W/V$. Armed with this terminology, we are ready to present the main result of this note.

\begin{theorem*}{\it Let $k$ be a field and suppose that $A$ is a $k$-algebra. Suppose further that $A$ does not have ideals of $k$-codimension one. Then no basis of $A$ is closed under multiplication.}\end{theorem*}

In the case where $A$ is a proper field extension of $k$ (which implies $\dim_k A\geqslant 2$), $A$ does not have ideals other than $\{0\}$ and $A$ (algebras with this property are called \emph{simple}), so in particular it does not have ideals of $k$-codimension 1. Indeed, if $a\in A$ is a non-zero element in $A$, then for any $b\in A$ one has $b=ba^{-1}a$, so $A$ itself is the smallest ideal containing $a$. Consequently, we derive Mabry's result as a corollary and the conclusion of our main result extends to the class of all simple $k$-algebras. It is well-known that the $k$-algebras $M_n(k)$ consisting of all $n\times n$-matrices over $k$ ($n\in \mathbb{N}$) are prototypical examples of simple $k$-algebras. Another familiar example is the $\mathbb{R}$-algebra of quaternions. On the other hand, the Cartesian products $k^n$ ($n\geqslant 1$) do have bases closed under multiplication, namely
$$(1, 0, 0, \ldots), (1, 1, 0, \ldots), \ldots, (1, 1, 1, \ldots). $$
Semigroup algebras constitute another class of $k$-algebras having bases closed under multiplication. Indeed, fix a field $k$ and a semigroup $(S,\cdot)$ (that is, a set with an associative binary operation). Denote by $k[S]$ the vector space over $k$ comprising all functions $f\colon S\to k$ that assume at most finitely many non-zero values. For each $s\in S$ define $\delta_s$ to be the function that takes value 1 on $s$ and 0 otherwise. Then the operation $\delta_s \cdot \delta_t = \delta_{st}$ ($s,t\in S$) extends uniquely to an associative operation on $k[S]$ that makes it a $k$-algebra. Clearly $\{\delta_s\colon s\in S\}$ is a basis for $k[S]$ that is closed under multiplication.\smallskip

\begin{proof}[Proof of the Main Theorem]Let $A$ be a $k$-algebra with a basis $H$ closed under multiplication. Consider the map $f\colon A\to k$ given by 
$$f(x)=\sum_{h\in H} c_h, $$
where $x=\sum_{h\in H}c_h h$ ($c_h\in k, h\in H$) is the unique basis expansion for $x\in A$. (This map is well-defined as there are at most finitely many non-zero scalars $c_h$ in the basis expansion of $x$.) Then $f$ is a linear functional on $A$. As $H$ is closed under multiplication, $f$ is multiplicative too, as already observed \cite[Problem 7 in Section 14.6]{kuczma} in the case where $A$ is the field of reals regarded as an algebra over the field of rational numbers. Indeed, given two elements $x=\sum_{h\in H}c_h h$ and $y\in \sum_{g\in H}d_g g$ in $A$, where $c_h, d_g\in k,\, h,g\in H$ (recall that there are at most finitely many non-zero scalars $c_h$ and $d_h$, where $h\in H$), we have 
$$f(x \cdot y ) = f( \sum_{h\in H }\sum_{g\in H}c_hd_g \underbrace{h\cdot g}_{\in H}) = \sum_{h\in H}\sum_{g\in H}c_hd_g = f(x)\cdot f(y), $$
as the basis expansion is unique.

In the case where $A$ has $k$-dimension at least two, there exist two distinct elements $h, h^\prime\in H$, so one has $f(h - h^\prime)=0$, hence $f$ has non-trivial kernel. As $f$ is non-zero, linear and multiplicative, it is a homomorphism of $k$-algebras. Consequently, the kernel of $f$ is a~proper two-sided ideal of $A$. Moreover, as $f$ is a surjection onto $k$, the $k$-codimension of $\ker f$ is equal to 1. \end{proof}\smallskip

\noindent \textbf{Acknowledgements.} We wish to thank the anonymous referee for valuable comments that improved the presentation of this note significantly. 

\end{document}